\documentclass[11pt]{article}
\usepackage{amsmath}
\usepackage{amssymb}
\usepackage{amsthm}
\usepackage{bookmark}
\usepackage{geometry}
\usepackage{hyperref}
\newtheorem{theorem}{Theorem}
\newtheorem{corollary}[theorem]{Corollary}
\newtheorem{lemma}[theorem]{Lemma}
\newtheorem{proposition}[theorem]{Proposition}
\theoremstyle{definition}
\newtheorem{example}[theorem]{Example}
\newtheorem{remark}[theorem]{Remark}
\DeclareMathOperator{\im}{im}
\DeclareMathOperator{\rad}{rad}
\renewcommand{\labelenumi}{\upshape(\roman{enumi})}
\renewcommand{\theenumi}\labelenumi

\title{Permutation polynomials from quadratic and linear forms}
\author{Ruikai Chen\\\small School of Mathematical Sciences, South China Normal University, Guangzhou 510631, China\\\small Email: \href{mailto:chen.rk@outlook.com}{chen.rk@outlook.com}; \href{mailto:rchen@scnu.edu.cn}{rchen@scnu.edu.cn}}
\date{}

\begin{document}

\maketitle

\begin{abstract}
We investigate a family of permutation polynomials of finite fields using the theory of quadratic and linear forms. For an odd prime power $q$, consider polynomials over $\mathbb F_{q^n}$ of the form $\sum_{0\le i\le j<n}a_{ij}x^{q^i+q^j}+L(x)$, where $a_{ij}\in\mathbb F_{q^n}$ and $L$ serves as an $\mathbb F_q$-linear endomorphism of $\mathbb F_{q^n}$. In this paper, we begin with a general treatment over arbitrary fields of characteristic not $2$, and then provide explicit characterizations for several classes of such permutation polynomials.

\textit{Keywords:} finite field, linear form, permutation polynomial, quadratic form
\end{abstract}

\section{Introduction}

Let $\mathbb F_q$ be a finite field and $n$ be a positive integer. Consider the extension $\mathbb F_{q^n}$, also viewed as a vector space over $\mathbb F_q$, and the polynomial ring $\mathbb F_{q^n}[x]$. For a polynomial over $\mathbb F_{q^n}$ where the exponent of each term is a power of $q$, we call it a $q$-linear polynomial over $\mathbb F_{q^n}$, since polynomials in this form induce all $\mathbb F_q$-linear endomorphisms of $\mathbb F_{q^n}$. A $q$-linear polynomial over $\mathbb F_{q^n}$ is exactly a polynomial $L$ such that $\mathrm{Tr}(uL(x))$ is a linear form of $\mathbb F_{q^n}/\mathbb F_q$ for every $u\in\mathbb F_{q^n}$, where $\mathrm{Tr}$ is the trace map of $\mathbb F_{q^n}/\mathbb F_q$. Similarly, another form of polynomials can be characterized as follows.

\begin{proposition}\label{polynomial}
For a polynomial $Q$ over $\mathbb F_{q^n}$ of degree less than $q^n$ with $q$ odd, $\mathrm{Tr}(uQ(x))$ is a quadratic form of $\mathbb F_{q^n}/\mathbb F_q$ for every $u\in\mathbb F_{q^n}$ if and only if
\[Q(x)=\sum_{0\le i\le j<n}a_{ij}x^{q^i+q^j}\]
for some $a_{ij}\in\mathbb F_{q^n}$.
\end{proposition}
\begin{proof}
The sufficiency is easy to verify, and the necessity follows from the results in Section \ref{injection}. Specifically, the trace form of $\mathbb F_{q^n}/\mathbb F_q$, defined by
\[(u,v)\mapsto\mathrm{Tr}(uv)\]
for $u,v\in\mathbb F_{q^n}$, is a nondegenerate symmetric bilinear form, and thus every quadratic form of $\mathbb F_{q^n}/\mathbb F_q$ can be expressed as $\mathrm{Tr}(x\ell(x))$ for some $q$-linear polynomial $\ell$ over $\mathbb F_{q^n}$. Now Proposition \ref{quadratic_map} gives the required form of the polynomial $Q$.
\end{proof}

Polynomials of this class, also known as Dembowski--Ostrom polynomials, play an important role in the theory of finite fields and their applications. When $q$ is even, these polynomials serve as permutation polynomials of $\mathbb F_{q^n}$ under certain conditions as in \cite{blokhuis2001permutations}. Further examples of related permutation polynomials can be found in \cite{tu2014several,ding2015permutation,gupta2016some,bhattacharya2017some,li2017new,zha2017further,wang2018six} and the references therein.

However, when $q$ is odd, the polynomial $Q$ is never a permutation polynomial of $\mathbb F_{q^n}$ since $Q(x)=Q(-x)$. To circumvent this, we may consider $Q+L$ for a $q$-linear polynomial $L$ over $\mathbb F_{q^n}$. Although some permutation polynomials of this form might arise from existing constructions, we aim to develop a general framework for studying their permutation properties through quadratic and linear forms. As a key to our approach, for every $u\in\mathbb F_{q^n}$, $\mathrm{Tr}(uQ(x))$ induces a quadratic form and $\mathrm{Tr}(uL(x))$ induces a linear form of $\mathbb F_{q^n}/\mathbb F_q$. The permutation
property can then be characterized in terms of the radicals of these quadratic forms and the associated linear forms.

In Section \ref{injection}, for the polynomials of interest, we abstract their essential features to study a class of maps on a bilinear space over an arbitrary field of characteristic not $2$, and establish a general criterion for them to be injective (see Theorem \ref{injective} below). In Section \ref{matrix}, we proceed to investigate those maps represented by some special matrices. Subsequently, in Section \ref{translation}, we introduce some results on polynomials associated with quadratic forms over finite fields and pave the way to Section \ref{permutation}, where several explicit families of permutation polynomials are derived directly in terms of their coefficients.

\section{A class of injections on bilinear spaces}\label{injection}

Let us first recall some fundamentals of the theory of bilinear spaces. Let $K$ be a field of characteristic different from $2$, and $V$ an $n$-dimensional vector space over $K$. Fix a nondegenerate symmetric bilinear form, denoted by $\langle\cdot,\cdot\rangle$. Now $V$ together with $\langle\cdot,\cdot\rangle$ forms a bilinear space.

Let $V^*$ be the dual space of $V$. There is a canonical isomorphism $\delta$ from $V$ to $V^*$ given by
\[v\mapsto\langle\cdot,v\rangle\]
for $v\in V$. For a basis $\beta_1,\dots,\beta_n$ of $V/K$, its dual basis $\beta_1^*,\dots,\beta_n^*$ is the basis of $V^*/K$ such that
\[v=\sum_{i=1}^n\beta_i^*(v)\beta_i,\]
for every $v\in V$.

For a $K$-linear endomorphism $\ell$ of $V$, there exists a unique $K$-linear endomorphism $\ell'$, called the adjoint of $\ell$ (with respect to the bilinear form), such that
\[\langle u,\ell(v)\rangle=\langle\ell'(u),v\rangle\]
for all $u,v\in V$. For a subspace $W$ of $V$, its orthogonal complement $W^\perp$ is the subspace consisting of $v\in V$ such that
\[\langle w,v\rangle=0\]
for all $w\in W$. It is easy to see that the adjoint relation and the orthogonal complement relation are both symmetric, with $(\ker\ell)^\perp=\im\ell'$.

Let $\mathcal Q$ be a quadratic form of $V/K$ and $\mathcal B$ be its associated symmetric bilinear form satisfying
\[\mathcal Q(u+v)-\mathcal Q(u)-\mathcal Q(v)=2\mathcal B(u,v)\]
for $u,v\in V$. This induces a $K$-linear map $\delta_\mathcal B$ from $V$ to $V^*$ that sends $v$ to $\mathcal B(\cdot,v)$. The radical of $\mathcal Q$ (or the radical of $\mathcal B$), denoted by $\rad\mathcal Q$ (or $\rad\mathcal B$), is defined as the kernel of $\delta_\mathcal B$. The image of $\delta_\mathcal B$ is precisely the subspace of those linear forms whose kernels contain $\rad\mathcal Q$, which is isomorphic to $(V/\rad\mathcal Q)^*$. Furthermore, $\mathcal B$ can be written as $\langle\cdot,\ell(\cdot)\rangle$ for a unique self-adjoint $K$-linear endomorphism $\ell$ of $V$, so that $\rad\mathcal Q=\ker\ell$.

Every quadratic form of $V/K$ can be diagonalized; i.e., there exists a basis $\beta_1,\dots,\beta_n$ of $V/K$ such that
\[\mathcal Q(v)=\sum_{i=1}^n\mathcal Q(\beta_i)\beta_i^*(v)^2,\]
and
\[\mathcal B(u,v)=\sum_{i=1}^n\mathcal Q(\beta_i)\beta_i^*(u)\beta_i^*(v),\]
for $u,v\in V$. Clearly, $\rad\mathcal Q$ is the subspace spanned by $\{\beta_i\mid\mathcal Q(\beta_i)=0\text{ for }1\le i\le n\}$. This also shows that every subspace of $V/K$ is the radical of some quadratic form. In particular, a quadratic form of $V/K$ is zero if and only if its radical is $V$.

\begin{theorem}\label{quadratic_linear}
For a quadratic form $\mathcal Q$ and a linear form $\mathcal L$ of $V/K$, take a basis $\beta_1,\dots,\beta_n$ of $V/K$ with $\mathcal Q(v)=\sum_{i=1}^n\mathcal Q(\beta_i)\beta_i^*(v)^2$ for $v\in V$. Then the following conditions are equivalent:
\begin{enumerate}
\item there is an index $i$ such that $\mathcal Q(\beta_i)=0$ and $\mathcal L(\beta_i)\ne0$;
\item $\rad\mathcal Q\not\subseteq\ker\mathcal L$;
\item $\delta^{-1}(\mathcal L)\notin(\rad\mathcal Q)^\perp$.
\end{enumerate}
\end{theorem}
\begin{proof}
Note that $\mathcal L(\rad\mathcal Q)$ is spanned by $\{\mathcal L(\beta_i)\mid\mathcal Q(\beta_i)=0\text{ for }1\le i\le n\}$. Hence, $\rad\mathcal Q\subseteq\ker\mathcal L$ if and only if $\mathcal L(\beta_i)=0$ for every index $i$ satisfying $\mathcal Q(\beta_i)=0$. This proves the equivalence between the first two conditions.

Next, we have $\big\langle v,\delta^{-1}(\mathcal L)\big\rangle=\mathcal L(v)$ by definition. Then $\rad\mathcal Q\subseteq\ker\mathcal L$ if and only if $\big\langle v,\delta^{-1}(\mathcal L)\big\rangle=0$ for every $v\in\rad\mathcal Q$; that is, $\delta^{-1}(\mathcal L)\in(\rad\mathcal Q)^\perp$.
\end{proof}

A map $L$ from $V$ to $V$ is a $K$-linear endomorphism if and only if $\lambda\circ L\in V^*$ for every $\lambda\in V^*$. This suggests a definition of ``$K$-quadratic'' maps on $V$.

\begin{proposition}\label{quadratic_map}
For a map $Q$ from $V$ to $V$ and a basis $\beta_1,\dots,\beta_n$ of $V/K$, the following conditions are equivalent:
\begin{enumerate}
\item\label{c1} $Q(kv)=k^2Q(v)$ for every $k\in K$ and $v\in V$, and the map from $V\times V$ to $V$ sending $(u,v)$ to $Q(u+v)-Q(u)-Q(v)$ is $K$-bilinear;
\item\label{c2} $\lambda\circ Q$ is a quadratic form of $V/K$ for every $\lambda\in V^*$;
\item\label{c3} $Q=\sum_{i=1}^n\beta_i\mathcal Q_i$ for some quadratic forms $\mathcal Q_1,\dots,\mathcal Q_n$ of $V/K$.
\end{enumerate}
\end{proposition}
\begin{proof}
One can easily verify that \ref{c1} implies \ref{c2} and that \ref{c3} implies \ref{c1}, given the definition of quadratic forms. For $v\in V$, expanding $Q(v)$ with respect to the basis yields
\[Q(v)=\sum_{i=1}^n(\beta_i^*\circ Q)(v)\beta_i,\]
so \ref{c2} implies \ref{c3}.
\end{proof}

From now on, let $Q$ be a map with these properties, and $L$ be a $K$-linear endomorphism of $V$.

\begin{theorem}\label{injective}
The map $Q+L$ is injective if and only if $\rad\lambda\circ Q\not\subseteq\ker\lambda\circ L$ for every nonzero $\lambda\in V^*$.
\end{theorem}
\begin{proof}
Define $D_u(v)=(Q+L)\big(\frac{u+v}2\big)-(Q+L)\big(\frac{u-v}2\big)$ and $B(u,v)=Q\big(\frac{u+v}2\big)-Q\big(\frac{u-v}2\big)$ for $u,v\in V$, such that $D_u$ is a $K$-linear endomorphism of $V$ and $B$ is a $K$-bilinear map on $V$; to see this, observe that $D_u(v)=B(u,v)+L(v)$ and
\[B(u,v)=\frac14(Q(u+v)-Q(u)-Q(v))-\frac14(Q(u-v)-Q(u)-Q(-v)).\]
Moreover, $\lambda\circ B$ is the associated symmetric bilinear form of $\lambda\circ Q$.

The following form a chain of equivalent statements:
\begin{itemize}
\item the map $Q+L$ is not injective;
\item $\ker D_u\ne\{0\}$ for some $u\in V$;
\item $\im D_u\subseteq\ker\lambda$ for some $u\in V$ and some nonzero $\lambda\in V^*$;
\item $-\lambda\circ L$ is equal to the linear form $(\lambda\circ B)(u,\cdot)$ for some $u\in V$ and some nonzero $\lambda\in V^*$;
\item $\rad\lambda\circ Q\subseteq\ker\lambda\circ L$ for some nonzero $\lambda\in V^*$.
\end{itemize}
The last equivalence holds because $\im\delta_{\lambda\circ B}$, as introduced earlier, consists precisely of those linear forms vanishing on $\rad\lambda\circ Q$.
\end{proof}

\begin{proposition}\label{bijective}
If the map $Q+L$ is injective, then $L$ is bijective.
\end{proposition}
\begin{proof}
If $L$ is not bijective, then $\im L\subseteq\ker\lambda$ for some nonzero $\lambda\in V^*$, for which $\ker\lambda\circ L=V$ and $\rad\lambda\circ Q\subseteq\ker\lambda\circ L$.
\end{proof}

For the map $Q+L$ to be injective, it is necessary that $L$ is a $K$-linear automorphism of $V$, in which case
\[Q+L=L\circ\big(L^{-1}\circ Q+\mathrm{Id}\big),\]
where $\mathrm{Id}$ is the identity map on $V$. Therefore, within the entire class of injections of the form $Q+L$, we may assume $L=\mathrm{Id}$ without loss of generality. Concerning only the quadratic part $Q$, there are several equivalent representations:
\begin{itemize}
\item a map $Q$ on $V$ such that $\lambda\circ Q$ is a quadratic form of $V/K$ for every $\lambda\in V^*$,
\item a tuple of $n$ quadratic forms of $V/K$, and
\item a $K$-linear map from $V^*$ to the space of quadratic forms of $V/K$ ($\lambda\mapsto\lambda\circ Q$).
\end{itemize}
Specifically, given quadratic forms $\lambda\circ Q$ with $\lambda$ ranging over a basis of $V^*/K$, it is straightforward to recover $Q$ using the dual basis.

Clearly, when $n=1$ the map $Q+L$ is injective only if $Q$ is zero. The following theorem asserts the existence of such an injection with nonzero $Q$ whenever $n>1$, and the proof provides a construction.

\begin{theorem}\label{existence}
If $n>1$, then there exists an injection on $V$ of the form $Q+\mathrm{Id}$ with $Q$ nonzero.
\end{theorem}
\begin{proof}
For an arbitrary subspace $W$ of $V$ with $0<\dim W<n$, the quadratic forms of $V/K$ whose radicals contain $W$ form a vector space $\mathcal F$ over $K$. Then there exists a nonzero $K$-linear map from $V^*$ to $\mathcal F$ whose kernel contains $\delta\big(W^\perp\big)$, since $\dim\mathcal F>0$ and $\dim W^\perp<n$. In other words, the map $Q$ can be chosen to be nonzero, such that $\rad\lambda\circ Q=V$ for $\lambda\in\delta\big(W^\perp\big)$ and $W\subseteq\rad\lambda\circ Q$ for every $\lambda\in V^*$. If $\lambda$ is nonzero and $\lambda\in\delta\big(W^\perp\big)$, then $\rad\lambda\circ Q=V\not\subseteq\ker\lambda$. If $\lambda\notin\delta\big(W^\perp\big)$, then $\lambda\notin\delta\big((\rad\lambda\circ Q)^\perp\big)$ since $(\rad\lambda\circ Q)^\perp\subseteq W^\perp$. Then $Q+\mathrm{Id}$ is injective by Theorem \ref{quadratic_linear} and Theorem \ref{injective}.
\end{proof}

\section{A matrix approach}\label{matrix}

In this section, we continue the discussion over the field $K$. Since all $n$-dimensional vector spaces over $K$ are isomorphic, we restrict our attention to $K^n$, working in the standard language of matrices. As is conventional, we write elements of $K^n$ as column vectors, linear forms of $K^n/K$ as row vectors, and $K$-linear endomorphisms of $K^n$ as matrices.

Given Proposition \ref{quadratic_map}, let
\begin{equation}\label{varphi}Q(v)=\begin{pmatrix}v^\mathrm TA_1v\\\vdots\\v^\mathrm TA_nv\end{pmatrix}\quad\text{and}\quad L(v)=Bv,\end{equation}
for $v\in K^n$, where $A_1,\dots,A_n$ are $n\times n$ symmetric matrices and $B$ is an $n\times n$ matrix over $K$. Then for every linear form $\lambda$ of $K^n/K$ written as $\begin{pmatrix}\lambda_1&\cdots&\lambda_n\end{pmatrix}$, the quadratic form $\lambda\circ Q$ is represented by the symmetric matrix $\sum_{i=1}^n\lambda_iA_i$, with $\rad\lambda\circ Q=\ker(\sum_{i=1}^n\lambda_iA_i)$. By Theorem \ref{injective}, the map $Q+L$ is an injection on $K^n$ if and only if, for every nonzero row vector $\lambda$ of dimension $n$, one has $\ker(\sum_{i=1}^n\lambda_iA_i)\not\subseteq\ker\lambda B$.

\begin{example}
Here is a concrete example on $K^2$. Let
\[A_1=A_2=\begin{pmatrix}1&1\\1&1\end{pmatrix}\quad\text{and}\quad B=\begin{pmatrix}1&0\\0&-1\end{pmatrix}\]
over $K$. Then for the row vector $\lambda=\begin{pmatrix}\lambda_1&\lambda_2\end{pmatrix}$ over $K$, we have
\[\lambda_1A_1+\lambda_2A_2=(\lambda_1+\lambda_2)\begin{pmatrix}1&1\\1&1\end{pmatrix},\]
and
\[\lambda B=\begin{pmatrix}\lambda_1&-\lambda_2\end{pmatrix}.\]
If $\lambda_1+\lambda_2=0$ with $\lambda_1\ne0$, then $\ker(\lambda_1A_1+\lambda_2A_2)=K^2$ and $\lambda B$ is nonzero. If $\lambda_1+\lambda_2\ne0$, then $\ker(\lambda_1A_1+\lambda_2A_2)$ is spanned by $\begin{pmatrix}1\\-1\end{pmatrix}$, while $\lambda B\begin{pmatrix}1\\-1\end{pmatrix}=\lambda_1+\lambda_2\ne0$. In either case, $\lambda B$ does not vanish on $\ker(\lambda_1A_1+\lambda_2A_2)$.
\end{example}

\begin{example}
For $1\le i\le n$, choose $A_i$ as an $n\times n$ symmetric matrix whose first $i$ rows and first $i$ columns are all zero, and let $B$ be the $n\times n$ identity matrix over $K$. For any nonzero row vector $\lambda$, let $j$ be the least index for which $\lambda_j\ne0$, and let $v$ be the column vector with $1$ as its $j$th entry and $0$ elsewhere. Then $v$ belongs to $\ker(\sum_{i=1}^n\lambda_iA_i)$ and, on the other hand, $\lambda Bv=\lambda_j\ne0$. Hence, $\lambda B$ does not vanish on $\ker(\sum_{i=1}^n\lambda_iA_i)$.
\end{example}

For more families of injections on $K^n$ involving certain special matrices, we begin with an auxiliary result.

Partition an $n\times n$ matrix $A$ over a field into $m\times m$ blocks with its rows and columns grouped according to the same index partition, where every diagonal block is a square matrix. The block matrix $A$ is treated block-wise, where $A_{ij}$ denotes the submatrix at the $i$th (block) row and $j$th (block) column. Furthermore, $A$ is said to be strictly upper triangular (at the block level) if all blocks on and below its diagonal are zero, and it is said to be strictly upper triangular up to block permutation similarity, if $P^{-1}AP$ is strictly upper triangular for some permutation matrix $P$ that permutes the block indices.

\begin{lemma}\label{block}
For the $m\times m$ block matrix $A$, write a row vector $\lambda$ of dimension $n$ in block form accordingly. Then $A$ is strictly upper triangular up to block permutation similarity if and only if
\begin{equation}\label{condition}\text{for every nonzero \(\lambda\), \((\lambda A)_j\) is zero but \(\lambda_j\) is nonzero for some index \(j\).}\end{equation}
\end{lemma}
\begin{proof}
Suppose that $A$ is indeed strictly upper triangular up to block permutation similarity. Rearrange the block indices so that $A$ is strictly upper triangular. For an arbitrary nonzero row vector $\lambda$, let $j$ be the least index such that $\lambda_j$ is nonzero. It is readily seen that $(\lambda A)_j$ is zero, as desired.

Next, we prove the sufficiency by induction on $m$. When $m=1$, this is obvious. For the inductive step, suppose that the sufficiency holds for any $(m-1)\times(m-1)$ block matrix whose diagonal blocks are square with $m>1$. Fix an arbitrary index $i$ and let $\lambda$ be nonzero with $\lambda_i$ being zero. The condition \eqref{condition} implies that $(\lambda A)_j$ is zero with $\lambda_j$ nonzero for some index $j$. For the principal block submatrix $A^{(i)}$ obtained by deleting the $i$th row and the $i$th column of $A$, and the row vector $\lambda^{(i)}$ obtained by deleting the $i$th column of $\lambda$, each block of $\lambda^{(i)}A^{(i)}$ corresponds to $(\lambda A)_j$ for some index $j$ different from $i$. This indicates that $A^{(i)}$ satisfies the condition \eqref{condition} in place of $A$, and, by the inductive hypothesis, $A^{(i)}$ is strictly upper triangular up to block permutation similarity. Therefore, there must be a zero row in $A^{(i)}$.

We proceed by showing that there is a zero row in $A$. Assume to the contrary that there is a nonzero block in every row of $A$. Construct a sequence of indices iteratively as follows.
\begin{itemize}
\item Start with an arbitrary index $s_1$.
\item At each step, given $s_i$, choose a zero row in $A^{(s_i)}$ and set $s_{i+1}$ to be its corresponding index in $A$.
\item Continue until some index occurs again.
\end{itemize}
This process terminates in finitely many steps. Let $I$ be the set of the indices after the first occurrence of the repeating index in the constructed sequence. Define a relation $\leftarrow$ on $I$: $i\leftarrow j$ if and only if $i$ is the successor of $j$ in the sequence. Clearly, this relation acts as a permutation of $I$ (for every $i\in I$, there exists a unique index $j\in I$ with $i\leftarrow j$, and vice versa). For $i,j\in I$, if $i\leftarrow j$, then, by assumption, the $i$th row of $A$ is nonzero, but it corresponds to a zero row in $A^{(j)}$, so $A_{ij}$ is the unique nonzero block in the $i$th row of $A$. Conversely, if $A_{ij}$ is nonzero for $i,j\in I$, then $i\leftarrow j$ by the uniqueness established above. There exists a nonzero row vector $\lambda$ such that $\lambda_iA_{ij}$ is nonzero for $i\in I$ with $i\leftarrow j$, and $\lambda_i$ is zero for every $i\notin I$. Now, $\lambda_j$ being nonzero implies $j\in I$ and $(\lambda A)_j=\lambda_iA_{ij}$, where $i$ is the unique index in $I$ such that $i\leftarrow j$. This leads to a contradiction.

Finally, take a zero row in $A$ and let $i$ be its index. As shown earlier, some permutation of the indices excluding $i$ makes $A^{(i)}$ strictly upper triangular. Afterwards, move the $i$th row and the $i$th column of $A$ to the last position, so that the resulting matrix becomes strictly upper triangular.
\end{proof}

For a diagonalized quadratic form, there is a simple criterion to determine whether a linear form vanishes on its radical, as shown in Theorem \ref{quadratic_linear}. This prompts us to examine the map $Q+L$ with $Q$ defined by $n$ diagonal matrices $A_1,\dots,A_n$ and $L$ defined by an arbitrary matrix $B$. Nevertheless, it is noteworthy that multiple quadratic forms in general need not be simultaneously diagonalizable.

\begin{theorem}\label{diagonal}
Suppose that $A_1,\dots,A_n$ are diagonal matrices over $K$, and let $D$ be the $n\times n$ matrix whose $(i,j)$ entry is the $j$th diagonal entry of $A_i$. Then the map $Q+L$ defined in \eqref{varphi} is injective if and only if $B$ is invertible and $B^{-1}D$ is strictly upper triangular up to permutation similarity.
\end{theorem}
\begin{proof}
By Proposition \ref{bijective}, $B$ is necessarily invertible. For a row vector $\lambda$ of dimension $n$, the matrix $\sum_{i=1}^n\lambda_iA_i$ is diagonal and its $j$th diagonal entry is $(\lambda D)_j$. The map $Q+L$ is injective if and only if, for every nonzero $\lambda$, there exists some index $j$ with $(\lambda D)_j=0$ and $(\lambda B)_j\ne0$. Substituting $\lambda$ for $\lambda B$ and applying Lemma \ref{block} entry-wise, we obtain the desired result.
\end{proof}

We turn to quadratic forms represented by circulant matrices. There is a $K$-algebra homomorphism $\varphi$ from the polynomial ring $K[x]$ onto the algebra of $n\times n$ circulant matrices over $K$, defined by
\[x\mapsto\begin{pmatrix}0&1&0&\cdots&0\\0&0&1&\cdots&0\\\vdots&\vdots&\vdots&\ddots&\vdots\\0&0&0&\cdots&1\\1&0&0&\cdots&0\end{pmatrix},\]
whose kernel is the ideal generated by $x^n-1$. In addition, let $\psi(f)$ denote the row vector of the coefficients (in ascending order) of $f$ reduced modulo $x^n-1$, the same as the first row of $\varphi(f)$.

\begin{lemma}
For polynomials $f$ and $g$ over $K$, $\ker\varphi(f)\subseteq\ker\psi(g)$ if and only if $\gcd(f,x^n-1)$ divides $g$.
\end{lemma}
\begin{proof}
Suppose that $\ker\varphi(f)\subseteq\ker\psi(g)$. Note that $f\cdot\frac{x^n-1}{\gcd(f,x^n-1)}\equiv0\mod x^n-1$, so
\[\im\varphi\bigg(\frac{x^n-1}{\gcd(f,x^n-1)}\bigg)\subseteq\ker\varphi(f)\subseteq\ker\psi(g).\]
Then $\psi(g)\varphi\Big(\frac{x^n-1}{\gcd(f,x^n-1)}\Big)$, as well as the first row of $\varphi(g)\varphi\Big(\frac{x^n-1}{\gcd(f,x^n-1)}\Big)$, is zero. Consequently, $g\cdot\frac{x^n-1}{\gcd(f,x^n-1)}\equiv0\mod x^n-1$ and $\gcd(f,x^n-1)$ divides $g$.

Conversely, if $\gcd(f,x^n-1)$ divides $g$, then $fh\equiv g\mod x^n-1$ for some polynomial $h$ over $K$; in this case $\ker\varphi(f)\subseteq\ker\varphi(g)$, and \textit{a fortiori} $\ker\varphi(f)\subseteq\ker\psi(g)$.
\end{proof}

Suppose that $x^n-1$ is separable over $K$ (the characteristic of $K$ does not divide $n$) with $x^n-1=\prod_{i=1}^mp_i$ for distinct irreducible polynomials $p_1,\dots,p_m$ over $K$. For each $i$, take polynomials $r_{i,1},\dots,r_{i,\deg p_i}$ that are linearly independent modulo $x^n-1$ over $K$ and all divisible by $\frac{x^n-1}{p_i}$; for example, by the Chinese remainder theorem, let
\[r_{i,j}\equiv x^{j-1}\mod p_i\quad\text{and}\quad r_{i,j}\equiv0\mod\frac{x^n-1}{p_i}\]
for $1\le j\le \deg p_i$. All these polynomials $r_{i,j}$ reduced modulo $x^n-1$ form a basis of $K[x]/(x^n-1)$ over $K$. Order the basis in blocks as follows:
\begin{equation}\label{basis}(r_{1,1},\dots,r_{1,\deg p_1}),\dots,(r_{m,1},\dots,r_{m,\deg p_m}).\end{equation}
Denote by $R$ the matrix formed by row vectors $\psi(r_{i,j})$ in this order. For a polynomial $f$ over $K$, the coefficient vector of $f$ modulo $x^n-1$ in the basis is $\psi(f)R^{-1}$ as $\psi(f)=\big(\psi(f)R^{-1}\big)R$. Partition $\psi(f)R^{-1}$ into blocks as above, so that $p_i$ divides $f$ if and only if the $i$th block of $\psi(f)R^{-1}$ is zero.

\begin{theorem}
Under the above setting, let $A_1,\dots,A_n$ be symmetric matrices with $A_1=\varphi(f_1),\dots,A_n=\varphi(f_n)$ for polynomials $f_1,\dots,f_n$ over $K$, and let $C$ be the $n\times n$ matrix whose $i$th row is given by $\psi(f_i)$. Now, all $n\times n$ matrices are viewed as $m\times m$ block matrices according to \eqref{basis}. Then the map $Q+L$ defined in \eqref{varphi} is injective if and only if $B$ is invertible and $RB^{-1}CR^{-1}$ is strictly upper triangular up to block permutation similarity.
\end{theorem}
\begin{proof}
For a nonzero row vector $\lambda$, the coefficient vector of $\sum_{i=1}^n\lambda_if_i$ reduced modulo $x^n-1$ in the basis \eqref{basis} is $\lambda CR^{-1}$. For a polynomial $g$ over $K$, it follows from the preceding lemma that $\ker(\sum_{i=1}^n\lambda_i\varphi(f_i))\not\subseteq\ker\psi(g)$ if and only if $\gcd(\sum_{i=1}^n\lambda_if_i,x^n-1)$ does not divide $g$; that is, by the previous discussion, $\big(\lambda CR^{-1}\big)_j$ is zero but $\big(\psi(g)R^{-1}\big)_j$ is nonzero for some block index $j$. Choosing $g$ such that $\psi(g)=\lambda B$, we complete the proof by an argument analogous to that in the proof of Theorem \ref{diagonal}.
\end{proof}

The assumption that $x^n-1$ be separable over $K$ is made for simplicity. When $x^n-1$ has repeated factors, determining polynomial divisibility modulo $x^n-1$ becomes more complicated. A refined approach with a more elaborate block structure would likely handle the inseparable case as well.

With these in mind, we now concentrate on finite fields, where injections correspond precisely to permutation polynomials. Henceforth, let $K=\mathbb F_q$ with $q$ odd, and denote by $\mathrm{Tr}$ the trace map of $\mathbb F_{q^n}/\mathbb F_q$.

Take a pair of bases $\alpha_1,\dots,\alpha_n$ and $\beta_1,\dots,\beta_n$ of $\mathbb F_{q^n}$ over $\mathbb F_q$ such that
\[v=\sum_{i=1}^n\mathrm{Tr}(\alpha_iv)\beta_i\]
for every $v\in\mathbb F_{q^n}$. Via an isomorphism between $\mathbb F_{q^n}$ and $\mathbb F_q^n$, we can construct permutation polynomials of $\mathbb F_{q^n}$ from permutations of $\mathbb F_q^n$, and vice versa. For $n\times n$ symmetric matrices $A_1,\dots,A_n$, let
\[\mathcal Q_i(x)=\begin{pmatrix}\mathrm{Tr}(\alpha_1x)&\cdots&\mathrm{Tr}(\alpha_nx)\end{pmatrix}A_i\begin{pmatrix}\mathrm{Tr}(\alpha_1x)\\\vdots\\\mathrm{Tr}(\alpha_nx)\end{pmatrix}\]
for $1\le i\le n$. Then the map $Q$ defined in \eqref{varphi} corresponds to the polynomial
\[\sum_{i=1}^n\mathcal Q_i(x)\beta_i.\]
Given a matrix as an $\mathbb F_q$-linear endomorphism of $\mathbb F_q^n$, its corresponding polynomial is derived in the same manner. To conclude, the matrices studied in this section actually induce permutation polynomials of $\mathbb F_{q^n}$.

\section{On polynomials inducing quadratic forms over finite fields}\label{translation}

We adapt the results of Section \ref{injection} into the context of finite fields. Regard $\mathbb F_{q^n}$ as a bilinear space over $\mathbb F_q$ equipped with the trace form of $\mathbb F_{q^n}/\mathbb F_q$. A $q$-linear polynomial $\ell$ over $\mathbb F_{q^n}$ is naturally identified with an $\mathbb F_q$-linear endomorphism of $\mathbb F_{q^n}$. In this sense, arithmetic in the polynomial ring $\mathbb F_{q^n}[x]$ is taken modulo $x^{q^n}-x$ implicitly, when only the induced maps on $\mathbb F_{q^n}$ are concerned.

Let
\[\ell(x)=\sum_{i=0}^{n-1}a_ix^{q^i},\]
with $a_i\in\mathbb F_{q^n}$. Its adjoint is
\[\ell'(x)=a_0x+\sum_{i=1}^{n-1}(a_ix)^{q^{n-i}}.\]
Define a matrix associated with $\ell$ as
\[M_\ell=\begin{pmatrix}a_0&a_1&\cdots&a_{n-1}\\a_{n-1}^q&a_0^q&\cdots&a_{n-2}^q\\\vdots&\vdots&\ddots&\vdots\\a_1^{q^{n-1}}&a_2^{q^{n-1}}&\cdots&a_0^{q^{n-1}}\end{pmatrix};\]
that is, an $n\times n$ matrix with $(i,j)$ entry $a_{j-i}^{q^{i-1}}$, where the index is taken modulo $n$ (see \cite{wu2013linearized} for more details on this matrix). The rank $r$ of $\ell$ is exactly the rank of $M_\ell$. A routine computation shows that $M_{\ell'}=M_\ell^\mathrm T$ and the entries of
\[\begin{pmatrix}b_0&\cdots&b_{n-1}\end{pmatrix}M_\ell,\]
for $b_0,\dots,b_{n-1}\in\mathbb F_{q^n}$, correspond to the coefficients of the $q$-linear polynomial $\sum_{i=0}^{n-1}b_i\ell(x)^{q^i}$ reduced modulo $x^{q^n}-x$. The leading $r\times r$ principal submatrix of $M_\ell$ is nonsingular. To see this, suppose that
\[\begin{pmatrix}b_0&\cdots&b_{r-1}&0&\cdots&0\end{pmatrix}M_\ell=\begin{pmatrix}0&\cdots&0&c_r&\cdots&c_{n-1}\end{pmatrix}\]
for some $c_r,\dots,c_{n-1}\in\mathbb F_{q^n}$. Then
\[\sum_{i=0}^{r-1}b_i\ell(x)^{q^i}=\sum_{i=r}^{n-1}c_ix^{q^i},\]
which must be zero since the kernel of $\sum_{i=r}^{n-1}c_i^{q^{n-r}}x^{q^{i-r}}$, of degree at most $q^{n-1-r}$, contains $\ker\ell$, of dimension $n-r$. Accordingly, the kernel of $\sum_{i=0}^{r-1}b_ix^{q^i}$ contains $\im\ell$, and thus $b_0,\dots,b_{r-1}$ are all zero.

The quadratic form of $\mathbb F_{q^n}/\mathbb F_q$ induced by $\mathrm{Tr}(x\ell(x))$ is equal to that of $\frac12\mathrm{Tr}(x(\ell(x)+\ell'(x)))$, where $\ell+\ell'$ is self-adjoint. Hence, its radical is $\ker(\ell+\ell')$. Every linear form of $\mathbb F_{q^n}/\mathbb F_q$ is induced by $\mathrm{Tr}(ux)$ for some $u\in\mathbb F_{q^n}$. Now, Theorem \ref{quadratic_linear} can be translated as follows.

\begin{theorem}\label{radical}
For a $q$-linear polynomial $\ell$ over $\mathbb F_{q^n}$ and an element $u$ of $\mathbb F_{q^n}$, the following conditions are equivalent:
\begin{enumerate}
\item the radical of the quadratic form $\mathrm{Tr}(x\ell(x))$ is not contained in the kernel of the linear form $\mathrm{Tr}(ux)$;
\item $u\notin\im(\ell+\ell')$;
\item the rank $r$ of $M_{\ell+\ell'}$ is less than $n$, and for the leading $(r+1)\times(r+1)$ principal submatrix $S$ of $M_{\ell+\ell'}$, the matrix obtained by substituting $\begin{pmatrix}u&u^q&\cdots&u^{q^r}\end{pmatrix}^\mathrm T$ for the first column of $S$ is nonsingular.
\end{enumerate}
\end{theorem}
\begin{proof}
The equivalence between the first two conditions is guaranteed by the equality $\ker(\ell+\ell')^\perp=\im(\ell+\ell')$. The case $r=n$ is trivial. Assume $r<n$ and consider $\sum_{i=0}^rc_ix^{q^i}$, where $c_i$ is the $(i+1,1)$ cofactor of $S$. This is a nonzero polynomial. Indeed, when $r=0$, we have $c_0=1$; when $0<r<n$, the cofactor $c_0$ is the $q$th power of the determinant of the leading $r\times r$ principal submatrix of $M_{\ell+\ell'}$. By definition, for any column vector $v$ over $\mathbb F_{q^n}$ of dimension $r+1$, $\begin{pmatrix}c_0&\cdots&c_r\end{pmatrix}v$ is the determinant of the matrix obtained by substituting $v$ for the first column of $S$. Any $r+1$ columns of $M_{\ell+\ell'}$ are linearly dependent over $\mathbb F_{q^n}$, so $\begin{pmatrix}c_0&\cdots&c_r&0&\cdots&0\end{pmatrix}M_{\ell+\ell'}$ is zero and $\sum_{i=0}^rc_i(\ell+\ell')(x)^{q^i}=0$. Comparing the dimensions, we conclude that $\im(\ell+\ell')$ is equal to the kernel of $\sum_{i=0}^rc_ix^{q^i}$, and the desired result follows.
\end{proof}

Given Proposition \ref{polynomial}, we focus on the polynomial
\[\sum_{i=0}^\nu L_i\Big(x^{q^i+1}\Big)+L(x),\]
where $\nu=\big\lfloor\frac n2\big\rfloor$ and each $L_i$, as well as $L$, is a $q$-linear polynomial over $\mathbb F_{q^n}$. For $u\in\mathbb F_{q^n}$, this gives rise to a quadratic form,
\[\mathrm{Tr}\Bigg(u\sum_{i=0}^\nu L_i\Big(x^{q^i+1}\Big)\Bigg)=\mathrm{Tr}\Bigg(\sum_{i=0}^\nu L_i'(u)x^{q^i+1}\Bigg),\]
and a linear form,
\[\mathrm{Tr}(uL(x))=\mathrm{Tr}(L'(u)x),\]
of $\mathbb F_{q^n}/\mathbb F_q$.

\begin{corollary}\label{criterion}
The polynomial $\sum_{i=0}^\nu L_i\Big(x^{q^i+1}\Big)+L(x)$ is a permutation polynomial of $\mathbb F_{q^n}$ if and only if the radical of $\mathrm{Tr}\Big(\sum_{i=0}^\nu L_i'(u)x^{q^i+1}\Big)$ is not contained in $\ker\mathrm{Tr}(L'(u)x)$ for every $u\in\mathbb F_{q^n}^\times$. If this is the case, then $L$ is an $\mathbb F_q$-linear automorphism of $\mathbb F_{q^n}$.
\end{corollary}

In what follows, we assume that the linear part is the identity map ($L(x)=x$) for simplicity. For the families of permutation polynomials examined below, one can verify that this causes no loss of generality.

\section{Permutation polynomials of explicit forms}\label{permutation}

\subsection{Polynomials with quadratic forms over intermediate fields}

Consider quadratic forms of $\mathbb F_{q^n}/\mathbb F_q$ that can also be defined over an intermediate field $\mathbb F_{q^m}$, where $m$ is a positive divisor of $n$. Denote by $\mathrm{Tr}_m$ the trace map of $\mathbb F_{q^n}/\mathbb F_{q^m}$.

\begin{lemma}\label{intermediate}
Let $\ell$ be a $q^m$-linear polynomial over $\mathbb F_{q^n}$. Then the radical of $\mathrm{Tr}_m(x\ell(x))$ as a quadratic form of $\mathbb F_{q^n}/\mathbb F_{q^m}$ is equal to the radical of $\mathrm{Tr}(x\ell(x))$ as a quadratic form of $\mathbb F_{q^n}/\mathbb F_q$. For $u\in\mathbb F_{q^n}$, the radical is contained in $\ker\mathrm{Tr}_m(ux)$ if and only if it is contained in $\ker\mathrm{Tr}(ux)$.
\end{lemma}
\begin{proof}
The adjoint of $\ell$ with respect to the trace form of $\mathbb F_{q^n}/\mathbb F_{q^m}$ is exactly its adjoint with respect to the trace form of $\mathbb F_{q^n}/\mathbb F_q$; that is, $\mathrm{Tr}_m(u\ell(x))=\mathrm{Tr}_m(\ell'(u)x)$ as well as $\mathrm{Tr}(u\ell(x))=\mathrm{Tr}(\ell'(u)x)$. As $\ell$ is arbitrary, this implies that, for any $\mathbb F_{q^m}$-linear subspace ($\im\ell$) of $\mathbb F_{q^n}$, it is contained in $\ker\mathrm{Tr}_m(ux)$ if and only if it is contained in $\ker\mathrm{Tr}(ux)$ ($\ell'(u)=0$). Furthermore, the radical of $\mathrm{Tr}_m(x\ell(x))$ is $\ker(\ell+\ell')$, the same as that of $\mathrm{Tr}(x\ell(x))$. 
\end{proof}

\begin{proposition}\label{2d}
Let $n=2m$. For $a,b,c\in\mathbb F_{q^n}$, the radical of $\mathrm{Tr}\big(ax^2+bx^{q^m+1}\big)$ is not contained in $\ker\mathrm{Tr}(cx)$ if and only if one of the following conditions holds:
\begin{itemize}
\item $a=\mathrm{Tr}_m(b)=0$ and $c\ne0$;
\item $4a^{q^m+1}=\mathrm{Tr}_m(b)^2\ne0$ and $2ac^{q^m}-\mathrm{Tr}_m(b)c\ne0$.
\end{itemize}
\end{proposition}
\begin{proof}
By Lemma \ref{intermediate}, this is reduced to the problem concerning quadratic forms of the $2$-dimensional space $\mathbb F_{q^n}/\mathbb F_{q^m}$. For $\ell(x)=ax+bx^{q^m}$ with $\ell(x)+\ell'(x)=2ax+\mathrm{Tr}_m(b)x^{q^m}$, we have
\[M_{\ell+\ell'}=\begin{pmatrix}2a&\mathrm{Tr}_m(b)\\\mathrm{Tr}_m(b)&2a^{q^m}\end{pmatrix}.\]
The required condition holds if and only if either the matrix $M_{\ell+\ell'}$ is zero and $c\ne0$, or $M_{\ell+\ell'}$ has rank $1$ and
\[\begin{vmatrix}c&\mathrm{Tr}_m(b)\\c^{q^m}&2a^{q^m}\end{vmatrix}\ne0,\]
by Theorem \ref{radical}.
\end{proof}

\begin{theorem}\label{2dp}
Let $n=2m$. For $q$-linear polynomials $L_0$ and $L_1$ over $\mathbb F_{q^n}$, the polynomial $L_0\big(x^2\big)+L_1\big(x^{q^m+1}\big)+x$ is a permutation polynomial of $\mathbb F_{q^n}$ if and only if
\begin{itemize}
\item $2\alpha L_0'(x)=\mathrm{Tr}_m(L_1'(x))$, and
\item $\ker\big(x^{q^m}-\alpha x\big)\subseteq\ker L_0'$
\end{itemize}
for some $\alpha\in\mathbb F_{q^n}$ with $\alpha^{q^m+1}=1$. Moreover, when $L_1'(x)=\sum_{i=0}^{m-1}a_ix^{q^i}$ for $a_i\in\mathbb F_{q^n}$, $\ker\big(x^{q^m}-\alpha x\big)\subseteq\ker\mathrm{Tr}_m\circ L_1'$ if and only if $a_i+\alpha^{q^i}a_i^{q^m}=0$ for $0\le i<m$.
\end{theorem}
\begin{proof}
By Corollary \ref{criterion} and Proposition \ref{2d}, the polynomial in question is a permutation polynomial of $\mathbb F_{q^n}$ if and only if, for every $u\in\mathbb F_{q^n}^\times$,
\begin{itemize}
\item $4L_0'(u)^{q^m+1}=\mathrm{Tr}_m(L_1'(u))^2$, and
\item $2L_0'(u)u^{q^m}-\mathrm{Tr}_m(L_1'(u))u\ne0$ whenever $L_0'(u)\ne0$.
\end{itemize}
Suppose that this is the case.

Write $2L_0'(x)=\sum_{i=0}^{n-1}a_ix^{q^i}$ and $\mathrm{Tr}_m(L_1'(x))=\sum_{i=0}^{n-1}b_ix^{q^i}$ for $a_i,b_i\in\mathbb F_{q^n}$ with $b_{i+m}=b_i^{q^m}$, where indices are taken modulo $n$. Since $4L_0'(u)^{q^m+1}=\mathrm{Tr}_m(L_1'(u))^2$ for every $u\in\mathbb F_{q^n}^\times$, we have the identity $4L_0'(x)^{q^m+1}=\mathrm{Tr}_m(L_1'(x))^2$ in terms of polynomials. Expanding this gives
\[\sum_{i=0}^{n-1}\sum_{j=0}^{n-1}a_ia_j^{q^m}x^{q^i+q^{j+m}}=\sum_{i=0}^{n-1}\sum_{j=0}^{n-1}b_ib_jx^{q^i+q^j};\]
that is,
\[\sum_{i=0}^{n-1}\sum_{j=0}^{n-1}a_ia_{j-m}^{q^m}x^{q^i+q^j}=\sum_{i=0}^{n-1}\sum_{j=0}^{n-1}b_ib_jx^{q^i+q^j},\]
by reducing the polynomial modulo $x^{q^n}-x$. Since $q$ is odd, the exponents $q^i+q^j$ for $0\le i\le j<n$ are distinct and less than $q^n$. This allows us to compare the coefficients, and it turns out that
\[a_ia_{i-m}^{q^m}=b_i^2,\]
and
\[a_ia_{j-m}^{q^m}+a_ja_{i-m}^{q^m}=2b_ib_j\]
when $i\ne j$. Then
\[2a_ia_jb_ib_j=a_ia_j\Big(a_ia_{j-m}^{q^m}+a_ja_{i-m}^{q^m}\Big)=a_i^2a_ja_{j-m}^{q^m}+a_j^2a_ia_{i-m}^{q^m}=a_i^2b_j^2+a_j^2b_i^2,\]
and
\[(a_ib_j-a_jb_i)^2=a_i^2b_j^2-2a_ia_jb_ib_j+a_j^2b_i^2=0.\]
If $L_0'(x)=0$, then necessarily $\mathrm{Tr}_m(L_1'(x))=0$. In any case, one has $2\alpha L_0'(x)=\mathrm{Tr}_m(L_1'(x))$ for some $\alpha\in\mathbb F_{q^n}$, and $\alpha^{q^m+1}=1$ since
\[4\alpha^{q^m+1}L_0'(x)^{q^m+1}=\mathrm{Tr}_m(L_1'(x))^{q^m+1}=\mathrm{Tr}_m(L_1'(x))^2.\]
As a result,
\[2L_0'(u)u^{q^m}-\mathrm{Tr}_m(L_1'(u))u=2L_0'(u)\big(u^{q^m}-\alpha u\big),\]
which is zero only if $L_0'(u)=0$. Thus, the kernel of $x^{q^m}-\alpha x$, which has $q^m$ elements, is contained in $\ker L_0'$. This proves the necessity, and the sufficiency can now be easily checked.

For the remaining part, note that
\[\mathrm{Tr}_m(L_1'(x))\equiv\sum_{i=0}^{m-1}a_ix^{q^i}+\sum_{i=0}^{m-1}a_i^{q^m}x^{q^{i+m}}\equiv\sum_{i=0}^{m-1}\Big(a_i+\alpha^{q^i}a_i^{q^m}\Big)x^{q^i}\mod x^{q^m}-\alpha x.\]
It follows that $\ker\big(x^{q^m}-\alpha x\big)\subseteq\ker\mathrm{Tr}_m\circ L_1'$ if and only if $\sum_{i=0}^{m-1}\Big(a_i+\alpha^{q^i}a_i^{q^m}\Big)x^{q^i}=0$, since $x^{q^m}-\alpha x$ is separable.
\end{proof}

\begin{corollary}
For $a,b,c\in\mathbb F_{q^2}$, the polynomial $ax^2+bx^{q+1}+cx^{2q}+x$ is a permutation polynomial of $\mathbb F_{q^2}$ if and only if $2a+b^q=b^{2q-1}+2c^q=0$.
\end{corollary}
\begin{proof}
Consider the $2$-dimensional space $\mathbb F_{q^2}/\mathbb F_q$ with $\mathrm{Tr}(x)=x+x^q$. From the previous discussion, the polynomial $ax^2+bx^{q+1}+cx^{2q}+x$ is a permutation polynomial of $\mathbb F_{q^2}$ if and only if
\begin{itemize}
\item $2\alpha(ax+c^qx^q)=\mathrm{Tr}(bx)$, and
\item $b+\alpha b^q=0$
\end{itemize}
for some $\alpha\in\mathbb F_{q^2}$ with $\alpha^{q+1}=1$. Assume $b\ne0$, given that the case $b=0$ is straightforward.

Suppose that such an element $\alpha$ exists. Then $\alpha b^{q-1}=-1$ and 
\[(2a+b^q)x+\big(b^{2q-1}+2c^q\big)x^q=b^{q-1}(\mathrm{Tr}(bx)-2\alpha(ax+c^qx^q))=0.\]
Conversely, if $2a+b^q=b^{2q-1}+2c^q=0$, then let $\alpha=-b^{1-q}$, so that
\[\mathrm{Tr}(bx)-2\alpha(ax+c^qx^q)=b^{1-q}\big((2a+b^q)x+\big(b^{2q-1}+2c^q\big)x^q\big)=0,\]
as claimed.
\end{proof}

\begin{remark}
By choosing a basis of $\mathbb F_{q^2}/\mathbb F_q$, the above polynomial can be written as a bivariate quadratic polynomial over $\mathbb F_q$. Permutation polynomials in this bivariate form have been classified in \cite[Proposition 3.4]{pang2025determination}. Theorem \ref{2dp} concerns arbitrary $q$-linear polynomials $L_0$ and $L_1$ over $\mathbb F_{q^n}$, and thereby applies in a more general setting, where the target polynomials correspond to quadratic polynomials in $n$ variables over $\mathbb F_q$.
\end{remark}

\begin{example}
Take an element $\alpha$ of $\mathbb F_{q^n}$ with $\alpha^{q^m+1}=1$, and let $L_0'(x)=\frac1{2\alpha}\mathrm{Tr}_m(L_1'(x))$ with
\[L_1'(x)=\sum_{i=0}^{m-1}a_i\Big(x^{q^i}-\alpha^{-q^i}x^{q^{i+m}}\Big)\]
for $a_i\in\mathbb F_{q^n}$. A routine check shows that the conditions in Theorem \ref{2dp} are satisfied for arbitrary $a_i$.
\end{example}

Theorem \ref{2dp} precisely describes permutation polynomials from quadratic forms of a $2$-dimensional space. We then proceed to the $3$-dimensional case.

\begin{proposition}\label{3d}
Let $n=3m$. For $a,b,c\in\mathbb F_{q^n}$, the radical of $\mathrm{Tr}\big(ax^2+bx^{q^m+1}\big)$ is not contained in $\ker\mathrm{Tr}(cx)$ if and only if one of the following conditions holds:
\begin{itemize}
\item $a=b=0$ and $c\ne0$;
\item $2a=b^{q^{2m}-q^m+1}\ne0$ and $2a^{q^m}c-bc^{q^m}\ne0$;
\item $2a\ne b^{q^{2m}-q^m+1}$, $4a^{q^{2m}+q^m+1}+b^{q^{2m}+q^m+1}-\mathrm{Tr}_m\big(ab^{2q^m}\big)=0$, and
\[\begin{vmatrix}c&b&b^{q^{2m}}\\c^{q^m}&2a^{q^m}&b^{q^m}\\c^{q^{2m}}&b^{q^m}&2a^{q^{2m}}\end{vmatrix}\ne0.\]
\end{itemize}
\end{proposition}
\begin{proof}
Consider $\ell(x)=ax+bx^{q^m}$ and the matrix
\[M_{\ell+\ell'}=\begin{pmatrix}2a&b&b^{q^{2m}}\\b&2a^{q^m}&b^{q^m}\\b^{q^{2m}}&b^{q^m}&2a^{q^{2m}}\end{pmatrix}\]
of rank $r$. As in the proof of Proposition \ref{2d}, it suffices to determine the rank of $M_{\ell+\ell'}$ and employ Theorem \ref{radical} directly.

We show that $r<2$ if and only if $2a=b^{q^{2m}-q^m+1}$. If $r<2$, then
\[4a^{q^m+1}-b^2=\begin{vmatrix}2a&b\\b&2a^{q^m}\end{vmatrix}=0,\]
and
\[\Big(b^{q^{2m}-q^m+1}-2a\Big)^{q^m}b^{q^{2m}}=b^{q^m+1}-2a^{q^m}b^{q^{2m}}=\begin{vmatrix}b&b^{q^{2m}}\\2a^{q^m}&b^{q^m}\end{vmatrix}=0;\]
in any case one has $2a=b^{q^{2m}-q^m+1}$. Conversely, if $2a=b^{q^{2m}-q^m+1}$, then either $a=b=0$ or $b\ne0$ and
\[\ell(x)+\ell'(x)=b^{q^{2m}-q^m+1}x+bx^{q^m}+b^{q^{2m}}x^{q^{2m}}=b^{-q^m}\mathrm{Tr}_m\Big(b^{q^{2m}+1}x\Big),\]
which implies $r<2$.

Also note that $\det M_{\ell+\ell'}=8a^{q^{2m}+q^m+1}+2b^{q^{2m}+q^m+1}-2\mathrm{Tr}_m\big(ab^{2q^m}\big)$. Then the proposition follows immediately.
\end{proof}

The conditions here are too complicated to allow a characterization for permutation polynomials as in Theorem \ref{2dp}. Instead, we focus on a special case where every quadratic form involved has rank either $0$ or $1$.

\begin{theorem}
Let $n=3m$. For $q$-linear polynomials $L_0$ and $L_1$ over $\mathbb F_{q^n}$, if $2L_0'(x)=L_1'(x)^{q^{2m}-q^m+1}$, then $2L_0'(x)=\alpha^{q^{2m}-q^m+1}\mathrm{Tr}_m(\ell(x))$ and $L_1'(x)=\alpha\mathrm{Tr}_m(\ell(x))$ for some $\alpha\in\mathbb F_{q^n}$ and $q$-linear polynomial $\ell$ over $\mathbb F_{q^n}$. In this case, the polynomial $L_0\big(x^2\big)+L_1\big(x^{q^m+1}\big)+x$ is a permutation polynomial of $\mathbb F_{q^n}$ if and only if either $\alpha=0$ or $\alpha\ne0$ and $\alpha^{-q^m}\mathbb F_{q^m}\subseteq\ker\mathrm{Tr}_m\circ\ell$.
\end{theorem}
\begin{proof}
Suppose that $2L_0'(x)=L_1'(x)^{q^{2m}-q^m+1}$, and neither $L_0'$ nor $L_1'$ vanishes on $\mathbb F_{q^n}$ (otherwise the result is trivial). Then the map $x^{q^{2m}-q^m+1}$ restricted to $\im L_1'$ is $\mathbb F_q$-linear, such that
\[(u+v)^{q^{2m}-q^m+1}=u^{q^{2m}-q^m+1}+v^{q^{2m}-q^m+1}\]
for all $u,v\in\im L_1'$. Raising the equation to the $(q^m+1)$th power, we get
\[(u+v)^2=u^2+u^{1-q^{2m}+q^m}v^{q^{2m}-q^m+1}+u^{q^{2m}-q^m+1}v^{1-q^{2m}+q^m}+v^2,\]
since $\big(q^{2m}-q^m+1\big)(q^m+1)\equiv2\mod q^n-1$. Let $w=u^{q^{2m}-q^m}v^{-q^{2m}+q^m}$ with $uv\ne0$, so that
\[2uv=(u+v)^2-u^2-v^2=uv\big(w^{-1}+w\big),\]
and thus $w=1$. Then $u^{q^m-1}=v^{q^m-1}$ and $u\in v\mathbb F_{q^m}$. This implies that $\im L_1'$ is contained in some $1$-dimensional $\mathbb F_{q^m}$-linear subspace of $\mathbb F_{q^n}$, and we may write $L_1'(x)=\alpha\mathrm{Tr}_m(\ell(x))$ as required.

Now we have $2L_0'(u)=\alpha^{q^{2m}-q^m+1}\mathrm{Tr}_m(\ell(u))$ with $\alpha\ne0$ for $u\in\mathbb F_{q^n}^\times$, so that
\[2L_0'(u)^{q^m}u-L_1'(u)u^{q^m}=\alpha^{1-q^{2m}}\mathrm{Tr}_m(\ell(u))\Big(\alpha^{q^m}u-\alpha^{q^{2m}}u^{q^m}\Big),\]
which is zero if and only if $\mathrm{Tr}_m(\ell(u))=0$ or $\alpha^{q^m}u\in\mathbb F_{q^m}$. By Proposition \ref{3d}, the polynomial $L_0\big(x^2\big)+L_1\big(x^{q^m+1}\big)+x$ is a permutation polynomial of $\mathbb F_{q^n}$ if and only if $2L_0'(u)^{q^m}u-L_1'(u)u^{q^m}\ne0$ whenever $L_0'(u)\ne0$; that is, $\mathrm{Tr}_m\circ\ell$ vanishes on $\alpha^{-q^m}\mathbb F_{q^m}$.
\end{proof}

When $\alpha\ne0$, the condition $\alpha^{-q^m}\mathbb F_{q^m}\subseteq\ker\mathrm{Tr}_m\circ\ell$ can be formulated as explicit equations concerning the coefficients of $\ell$, as in Theorem \ref{2dp}.

\begin{example}
Take a nonzero element $\alpha$ of $\mathbb F_{q^n}$, and let
\[\ell(x)=\sum_{i=0}^{m-1}a_i\Big(x^{q^i}-\alpha^{q^{i+2m}-q^{i+m}}x^{q^{i+m}}\Big)\]
for $a_i\in\mathbb F_{q^n}$, so that $\ell$ vanishes on $\alpha^{-q^m}\mathbb F_{q^m}$. Then \textit{a fortiori} $\alpha^{-q^m}\mathbb F_{q^m}\subseteq\ker\mathrm{Tr}_m\circ\ell$.
\end{example}

\subsection{Polynomials involving trace maps}

Finally, we provide further constructions of permutation polynomials using the trace maps. The following result illustrates the constructive proof of Theorem \ref{existence}, where the radicals of all the quadratic forms involved have a nontrivial intersection.

\begin{theorem}
For a $q$-linear polynomial $L$ over $\mathbb F_{q^n}$ and $b\in\mathbb F_{q^n}^\times$, the polynomial $b\mathrm{Tr}(xL(x))+x$ is a permutation polynomial of $\mathbb F_{q^n}$ if and only if $(L+L')(b)=0$.
\end{theorem}
\begin{proof}
For $u\in\mathbb F_{q^n}$, note that
\[\mathrm{Tr}(ub\mathrm{Tr}(xL(x)))=\mathrm{Tr}(bu)\mathrm{Tr}(xL(x)),\]
and the radical of the quadratic form is $\ker(L+L')$ when $\mathrm{Tr}(bu)\ne0$. As in Corollary \ref{criterion}, the polynomial under consideration is a permutation polynomial of $\mathbb F_{q^n}$ if and only if $u\notin\ker(L+L')^\perp$ for every $u\in\mathbb F_{q^n}$ such that $\mathrm{Tr}(bu)\ne0$; i.e., $\im(L+L')\subseteq b^{-1}\ker\mathrm{Tr}$. Indeed, one has $\mathrm{Tr}(b(L(x)+L'(x)))=\mathrm{Tr}((L(b)+L'(b))x)$, as desired.
\end{proof}

\begin{remark}
The sufficient condition in the preceding theorem can be explained using the theory of \textit{linear translators} as in \cite{charpin2009does} or \cite{kyureghyan2011constructing}. If $(L+L')(b)=0$, then $b$ belongs to the radical of the quadratic form $\mathrm{Tr}(xL(x))$ of $\mathbb F_{q^n}/\mathbb F_q$; in other words, $b$ acts as a $0$-linear translator for the polynomial $\mathrm{Tr}(xL(x))$, giving rise to a permutation polynomial in the above form (see, e.g., \cite[Theorem 8]{kyureghyan2011constructing}). In contrast, our result provides a complete characterization for such a class of polynomials related to quadratic forms.
\end{remark}

\begin{example}
It suffices to find a $q$-linear polynomial $L$ over $\mathbb F_{q^n}$ with $\ker(L+L')\ne\{0\}$. To be specific, let $q=3$, $n=3$, and $\theta$ be a primitive element of $\mathbb F_{27}$ satisfying $\theta^3=\theta+1$. For $L(x)=\theta x+(2\theta+1)x^3$, we have
\[L(x)+L'(x)=2\theta x+(2\theta+1)x^3+(2\theta+1)^9x^9=2\theta x+(2\theta+1)x^3+(2\theta+2)x^9,\]
so that $(L+L')(1)=0$. Then $\mathrm{Tr}\big(\theta x^2+(2\theta+1)x^4\big)+x$ is a permutation polynomial of $\mathbb F_{27}$.
\end{example}

\begin{theorem}
For a $q$-linear polynomial $L$ over $\mathbb F_{q^n}$ with $a,b\in\mathbb F_{q^n}^\times$, the polynomial $\mathrm{Tr}(ax)L(bx)+x$ is a permutation polynomial of $\mathbb F_{q^n}$ if and only if $L'(a)=0$ and $\ker(bL'(x)-cx)=\{0\}$ for every $c\in\mathbb F_q^\times$.
\end{theorem}
\begin{proof}
For $u\in\mathbb F_{q^n}^\times$, we have
\[\mathrm{Tr}\big(u\mathrm{Tr}(ax)L(bx)\big)=\mathrm{Tr}(ax)\mathrm{Tr}(bL'(u)x).\]
When $L'(u)=0$, this quadratic form is zero. When $L'(u)\ne0$, choosing a basis adapted to the two nonzero linear forms $\mathrm{Tr}(ax)$ and $\mathrm{Tr}(bL'(u)x)$, according as the linear forms are linearly dependent or independent, we conclude that the intersection of their kernels is precisely the radical of the quadratic form. Moreover, one has
\[\ker\mathrm{Tr}(ax)\cap\ker\mathrm{Tr}(bL'(u)x)=(a\mathbb F_q+bL'(u)\mathbb F_q)^\perp.\]
Hence, the polynomial $\mathrm{Tr}(ax)L(bx)+x$ is a permutation polynomial of $\mathbb F_{q^n}$ if and only if $u\notin a\mathbb F_q+bL'(u)\mathbb F_q$ for every $u\in\mathbb F_{q^n}\setminus\ker L'$; i.e., $u\in a\mathbb F_q+bL'(u)\mathbb F_q$ implies $u\in\ker L'$ for every $u\in\mathbb F_{q^n}^\times$. 

Suppose that this is the case. Then clearly $a\in\ker L'$. If $u\in\ker(bL'(x)-cx)$ for $c\in\mathbb F_q^\times$, then $cu=bL'(u)$ and $u\in bL'(u)\mathbb F_q$, which implies $u\in\ker L'$ and $u=c^{-1}bL'(u)=0$. As $\ker(bL'(x)-cx)=\{0\}$, this proves the necessity.

Conversely, suppose that $L'(a)=0$ and $\ker(bL'(x)-cx)=\{0\}$ for every $c\in\mathbb F_q^\times$. If $u\in a\mathbb F_q+bL'(u)\mathbb F_q$, then either $u\in a\mathbb F_q$ or $cu-bL'(u)\in a\mathbb F_q$ for some $c\in\mathbb F_q^\times$. In the latter case,
\[cL'(u)-L'(bL'(u))=L'(cu-bL'(u))=0,\]
and thus $bL'(u)\in\ker(bL'(x)-cx)$. This means $L'(u)=0$ whenever $u\in a\mathbb F_q+bL'(u)\mathbb F_q$. The result then follows.
\end{proof}

\begin{remark}
Let $f(x)=\mathrm{Tr}(ax)L(bx)+x$ and $T(x)=\mathrm{Tr}(ax)$. Then
\[(T\circ f)(x)=T(x)\mathrm{Tr}(L'(a)bx)+T(x).\]
If $L'(a)=0$, then $T\circ f=T$, and, for each $c\in\mathbb F_q$, $f$ acts on $T^{-1}(c)$ as $cL(bx)+x$. Therefore, applying \cite[Lemma 1.2]{akbary2011constructing} we may reduce the permutation property of $f$ to the injectivity of these linear maps, which yields the same sufficient condition described above. Our theorem further establishes the necessity of the condition $L'(a)=0$.
\end{remark}

\begin{example}
Let $L'(x)=\mathrm{Tr}(\beta x)$ for some $\beta\in\mathbb F_{q^n}$, and choose $a,b\in\mathbb F_{q^n}^\times$ such that $\mathrm{Tr}(\beta a)=\mathrm{Tr}(\beta b)=0$. If $bL'(u)-cu=0$ for some $c\in\mathbb F_q^\times$ and $u\in\mathbb F_{q^n}$, then $b^{-1}u=c^{-1}\mathrm{Tr}(\beta u)\in\mathbb F_q$, and further,
\[b^{-1}cu=\mathrm{Tr}\big(\beta bb^{-1}u\big)=\mathrm{Tr}(\beta b)b^{-1}u=0.\]
Therefore, $L'(a)=0$ and $\ker(bL'(x)-cx)=\{0\}$ for every $c\in\mathbb F_q^\times$.
\end{example}

\begin{example}
Let $L'(x)=x^q-a^{q-1}x$ with $a,b\in\mathbb F_{q^n}^\times$ satisfying $a^{1-q}b^{-1}\in\mathbb F_q$, and suppose that $\gcd(n,q-1)=1$. Then $bL'(x)-cx=bx^q-\big(ba^{q-1}+c\big)x$ for $c\in\mathbb F_q^\times$, where
\[\big(a^{q-1}+b^{-1}c\big)^\frac{q^n-1}{q-1}=\big(1+a^{1-q}b^{-1}c\big)^\frac{q^n-1}{q-1}=\big(1+a^{1-q}b^{-1}c\big)^n\ne1,\]
since the $n$th power map is bijective on $\mathbb F_q^\times$. This implies $L'(a)=0$ and $\ker(bL'(x)-cx)=\{0\}$.
\end{example}

\section*{Acknowledgments}

The author thanks the anonymous reviewers for their comments, which have significantly improved the quality of this paper.

\bibliographystyle{abbrv}
\bibliography{references}

\end{document}